\input amstex.tex
\input amsppt.sty
\input epsf
\documentstyle{amsppt}
\magnification=\magstep{1.5} 
\pagewidth{6.5truein} \pageheight{9.5truein} \TagsOnRight
\NoBlackBoxes

\let\phi\varphi
\define\N{{\Bbb N}}
\define\R{{\Bbb R}}

\define\FF{{\Cal F}}

\define\CC{{\Cal C}}

\def\ss{\subset}

\def\a{\alpha}

\topmatter
\title Ernest Michael and theory of continuous selections \endtitle
\author Du\v{s}an Repov\v{s} and Pavel V. Semenov \endauthor

\leftheadtext{Du\v{s}an Repov\v{s} and Pavel V. Semenov}
\rightheadtext{Ernest Michael and theory of selections}

\address
Institute of Mathematics, Physics and Mechanics, and Faculty of Education, University of
Ljubljana, Jadranska 19, P. O. Box 2964, Ljubljana, Slovenia 1001
\endaddress
\email dusan.repovs\@guest.arnes.si
\endemail
\address
Department of Mathematics, Moscow City Pedagogical University,
2-nd Selsko\-khozyast\-vennyi pr.\,4,\,Moscow,\,Russia 129226
\endaddress
\email pavels\@orc.ru
\endemail


\endtopmatter

\document
\rightline{\it To follow the thoughts of a great man is the
most interesting science.} \rightline{\it A. S. Pushkin}

\bigskip
\bigskip

\head {\bf 1. Introduction} \endhead For a large number of those
working in topology, functional analysis, multivalued analysis,
approximation theory, convex geometry, mathematical economics,
control theory, and several  other areas, the year 1956
has always been strongly connected with the publication by Ernest
Michael of two fundamental papers on continuous selections which
appeared in the Annals of Mathematics \cite{4} \cite{5}.

With sufficient precision that
year marked the beginning
of the theory of continuous selections of multivalued mappings. In
the last fifty years the approach to multivalued mappings and
their selections, set forth by Michael \cite{4} \cite{5},
has well established itself in contemporary mathematics. Moreover,
it has become an indispensable tool for many mathematicians
working in vastly different areas.

Clearly, the principal reason for this is the naturality of the
concept of selection. In fact, many mathematical assertions can be
reduced to using
the linguistic reversal
``$\forall x \in X \quad \exists y \in Y \ldots$''.
However, as
soon as we speak of the
validity of  assertions of the
type
$$
\forall x \in X \quad \exists y \in Y \quad P(x,y)
$$
it is natural to associate to every $x$ a nonempty set of all
those $y$ for which $P(x,y)$ is true. In this way we obtain a
multivalued map which can be interpreted as a mapping, which
associates to every initial data $x \in X$ of some problem $P$ a
nonempty set of solutions of this problem
$$
F: x \mapsto \{y \in Y : P(x,y)\}, \quad F: X \to Y.
$$
The question of the
existence of selections in such a
setting turns out
to be the question about the
unique choice of
the
solution of the problem
under given initial conditions. Different types of selections are
considered in different mathematical categories.

One could say that the key importance of Michael's theory is not
so much in providing a comprehensive solution of diverse selection
problems in the category of topological spaces and continuous
maps, but rather
the immediate inclusion of the
obtained results into
the general context of development of topology. In a remarkable
number of cases, results of Michael on solvability of
the
selection
problems turned out to be the final answers, i.e. they provided
conditions which turned out to be necessary and sufficient.

Initially we were planning to write a survey paper, which would present the
development of the theory in the last half of the
 century and its many applications.
However, already our first attempts at such a project showed that the volume of
such a survey would invariably fill an entire book, hence it
would be inappropriate for this special issue.

After some deliberation we decided to limit ourselves to a survey of
only the papers of Michael on the
theory of selections and their
mutual
relations. For analogous reasons we do not give any precise
references to many developments in the theory of selections -- the
number of papers in this area is by now around one thousand. A considerable
number of facts on selections and theorems, which go beyond the
present paper, can be found in books and surveys \cite{R1-R15}
listed at the end of the paper.

\head {\bf 2. Bibliography} \endhead

Papers in scientific journals usually end with the list of
references. In our opinion, is it most reasonable to begin a
survey dedicated to the work of a single person, on one special
topic, spanning over 50 years, with a complete list of his papers
on the subject.

\bigskip
\centerline{\smc List of all papers by E. Michael on selections}
\bigskip

{\eightpoint
\roster

\item"1." {\it Topologies on spaces of subsets},
Trans. Amer. Math. Soc. {\bf 71} (1951), 152--182.

\item"2." {\it Selection theorems for continuous functions},
Proc. Int. Congr. Math. {\bf 2} (1954), 241--242.

\item"3." {\it Selected selection theorems},
Amer. Math. Monthly {\bf 63} (1956), 233--238.

\item"4." {\it Continuous selections I}, 
Ann. of Math. (2) {\bf 63} (1956), 361--382.

\item"5." {\it Continuous selections II}, 
Ann. of Math. (2) {\bf 64} (1956), 562--580.

\item"6." {\it Continuous selections III}, 
Ann. of Math. (2) {\bf 65} (1957), 375--390.

\item"7." {\it A theorem on semi-continuous set valued functions}, 
Duke Math. J. {\bf 26} (1959), 647--652.

\item"8." {\it Dense families of continuous selections}, 
Fund. Math. {\bf 47} (1959), 174--178.

\item"9." {\it Paraconvex sets}, 
Math. Scand. {\bf 7} (1959), 372--376.

\item"10." {\it Convex structures and continuous selections}, 
Canadian J. Math. {\bf 11} (1959), 556--575.

\item"11." {\it Continuous selections in Banach spaces}, 
Studia Math. Ser. Spec. (1963), 75--76.

\item"12." {\it A linear mapping between function spaces}, 
Proc. Amer. Math. Soc. {\bf 15} (1964), 407--409.

\item"13." {\it Three mapping theorems}, 
Proc. Amer. Math. Soc. {\bf 15} (1964), 410--415.

\item"14." {\it A short proof of the Arens-Eells embedding theorem}, 
Proc. Amer. Math. Soc. {\bf 14} (1964), 415--416.

\item"15." {\it A selection theorem}, 
Proc. Amer. Math. Soc. {\bf 17} (1966), 1404--1406.

\item"16." {\it Topological well-ordering},
Invent. Math. {\bf 6} (1968), 150--158.
(with R. Engelking and R. Heath)

\item"17." {\it A unified theorem on continuous selections}, 
Pacific J. Math. {\bf 87} (1980), 187--188.
(with C. Pixley)

\item"18." {\it Continuous selections and finite-dimensional sets}, 
Pacific J. Math. {\bf 87} (1980), 189--197.

\item"19." {\it Continuous selections and countable sets}, 
Fund. Math. {\bf 111} (1981), 1--10.

\item"20." {\it A parametrization theorem}, 
Topology Appl. {\bf 21} (1985), 87--94.
(with G. M\"agerl and R. D. Mauldin)

\item"21." {\it A note on a selection theorem}, 
Proc. Amer. Math. Soc. {\bf 99} (1987), 575--576.

\item"22." {\it Continuous selections avoiding a set}, 
Topology Appl. {\bf 28} (1988), 195--213.

\item"23." {\it A generalization of a theorem on continuous selections},
Proc. Amer. Math. Soc. {\bf 105} (1989), 236--243.

\item"24." {\it Some problems},  
Open Problems in Topology, J.\ van Mill and G.\ M.\ Reed, Editors,  
North--Holland, Amsterdam, 1990, pp.\ 271--277.

\item"25." {\it Some refinements of a selection theorem with $0$-dimensional domain}, 
Fund. Math. {\bf 140} (1992), 279--287.

\item"26." {\it Selection theorems with and without dimensional restriction}, 
Recent Developments of General Topology and its Applications, 
International Conference in Memory of Felix Hausdorff (1868--1942),  
Math. Res.  67, Berlin, 1992.

\item"27." {\it Representing spaces as images of $0$-dimensional spaces}, 
Topology Appl. {\bf 49} (1993), 217--220.
(with M. M. Choban)

\item"28." {\it A note on global and local selections}, 
Topology Proc. {\bf 18} (1993), 189--194.

\item"29." {\it A theorem of Nepomnyashchii on continuous subset-selections},
Topology Appl. {\bf 142} (2004), 235--244.

\item"30." {\it Continuous Selections}, 
Encyclopedia of General Topology, c--8 (2004), 107--109.

\endroster
}

We have selected the
papers on
selections
\cite{4}, \cite{5} and \cite{7}  to serve as the basis of the classification of the
entire list. Here is a reasonably precise diagram of relationship
among the papers from the list:

\centerline{
\epsfxsize=12cm
\epsffile{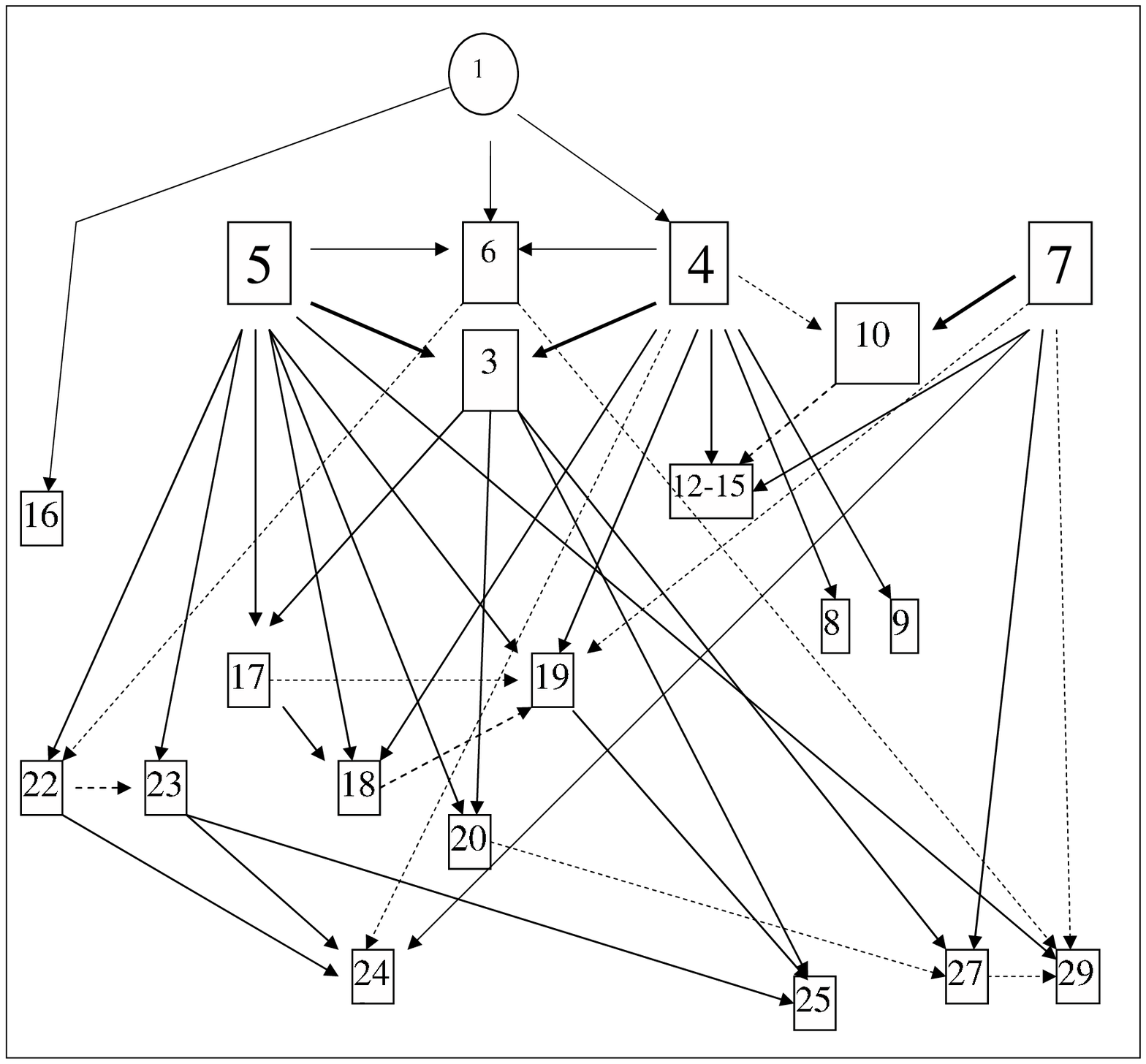}
}

\vfill\eject

Here the usual arrow means 
direct correlation and the dotted
arrow means an
implicit one.

Papers \cite{2}, \cite{11}, \cite{26}, \cite{28}, \cite{30} are not
included in this diagram, since they are either short announcements
(or abstracts) on conferences or they are devoted to popularization
of the subject.

\head {\bf 3. Papers from 1956} \endhead

A considerable number of fundamental mathematical papers can
be divided into two types. In such papers, as a rule, a
significant new theory is constructed or an important problem is
solved. This division is of course, conditional -- on the one
hand, in constructions of new theories one often encounters
difficult problems, on the other hand a solution of a difficult
problem often gives rise to a development of a significant new
theory.

The papers \cite{4} and \cite{5} are a clear cut
example of such a division. In \cite{4} an essentially new mathematical
theory is constructed, in the form of a branched tree, which unifies a
large number of sufficiently different theorems.
To the
contrary, in \cite{5} the principal result consists of the proof of a single
highly nontrivial theorem and all assertions and constructions
in this paper are devoted to the solution of this problem.

Another, linguistic difference between \cite{4} and \cite{5} is
connected with the notion of convexity: the formulations of
practically all theorems of \cite{4} use the term
{\sl convex}, where to the contrary, this word is practically
absent from \cite{5}. Finally, in \cite{4} Lebesgue
dimension is never used, while in \cite{5}, there are dimension
restrictions on the domains of the multivalued maps
everywhere.

One can say, with sufficient accuracy, that in \cite{5} the
finite-dimensional, purely topological analogue of such a
nontopological notions as convexity and local convexity are presented and studied. Without
any doubt, the best known assertion of \cite{4} is Theorem $3.2''$.

\proclaim{Theorem 1} The following properties of a  $T_{1}$-space
$X$ are equivalent:
\itemitem{(a)} $X$ is paracompact; and
\itemitem{(b)} If $Y$ is a Banach space, then every lower
semicontinuous (LSC) carrier $\phi:X \to \FF_{c}(Y)$ admits a
singlevalued continuous selection.
\endproclaim

Here $\FF_{c}(Y)$ denotes the family of all nonempty closed convex
subsets of $Y$. Observe that in \cite{4} Michael originally used
the term "carrier" instead of "multivalued mapping". In
mathematical practice the implication $(a) \Rightarrow (b)$ has the
widest application and is in folklore known as the
"Convex-Valued
Selection Theorem". The implication $(b) \Rightarrow (a)$ gives a
selection characterization of paracompactness.

The unusual numeration $3.2''$ for the theorem has a very
simple
explanation. In chapter 3 of \cite{4} Michael started by a
citation of Theorem 3.1 (Urysohn, Dugundji, Hanner) and Theorem
3.2 (Dowker) on the extensions of {\it singlevalued} mappings and
then presented the sequences:
\medskip
Theorem \,\,3.1,  Theorem \,\,3.1', Theorem
\,\,3.1'', Theorem \,\,3.1''' 
\medskip
Theorem \,\,3.2, Theorem \,\,3.2', Theorem
\,\,3.2'' 
\medskip
\noindent
of their analogs for {\it multivalued} mappings. To be more clear,
let us unify  Theorems 3.1 (a,b,c below) and $3.1'$ (a,d,e below) as
follows:

\proclaim{Theorem 2} The following properties of a $T_{1}$-space
$X$ are equivalent:
\itemitem{(a)} $X$ is normal;
\itemitem{(b)} The real line $\R$ is an extensor for $X$;
\itemitem{(c)} Every separable Banach space is extensor for $X$;
\itemitem{(d)} Every LSC carrier $\phi:X \to \CC(\R)$ admits a
singlevalued continuous selection; and
\itemitem{(e)} If $Y$ is a separable Banach space, then every LSC
carrier $\phi:X \to \CC(Y)$ admits a singlevalued continuous
selection.
\endproclaim

Thus, playing with words, such a series shows that the selection
theory in fact,
extends the theory of extensors. Here $\CC(Y) = \{Z
\in \FF_{c}(Y): Z$ is compact or $Z=Y\}$.

As asserted by Michael, Theorem $3.2''$ was his very first selection theorem,
the initial goal of which were generalizations of a theorem due to R. Bartle
and L. Graves on sections of linear continuous surjections between Banach spaces.
In
particular, Proposition 7.2 of \cite{4} states that such a section
can be chosen in an "almost" linear fashion (scalar homogenous) and
with the pointwise norm arbitrarily close
to the "minimal" of all
possible.

Thus the remaining Theorems $3.1''$--$3.2'$ are selection
characterizations of other properties of the domain of a
convex-valued mapping: normality, collectionwise normality,
normality and countable paracompactness, and perfect normality.
Many constructions and ideas from \cite{4} later became the basis
for subsequent research. For example, Lemma 5.2 in \cite{4} was
the first result in finding
pointwise dense families of
selections.

In comparison with \cite{4}, the paper \cite{5} originally dealt
only with the unique Theorem 1.2, the so-called "Finite-dimensional
selection theorem":

\proclaim{Theorem 3} Let $X$ be a paracompact space, $A\subset X$
a closed subset with $\dim_{X}(X\setminus A) \leq n+1$, $Y$ a
complete metric space, $\FF$ an equi-$LC^{n}$ family of nonempty
closed subsets of $Y$ and $\phi:X \to \FF$ an LSC map. Then every
singlevalued continuous selection of $\phi|_{A}$ can be extended to 
a singlevalued continuous selection of $\phi|_{U}$, for some
open subset $U\supset A$. If additionally every member of $\FF$ is
$n$-connected (briefly, $C^{n}$) then one can take $U=X$.
\endproclaim

Without any doubt, this is one of the most
complicated topological theorems, the six--step proof in
\cite{5} is  clearly a mathematical masterpiece. Various efforts
were made by several people in the last 50 years to simplify
this proof (or ``improve'' it), including ourselves.
However,
none of these versions turned out to be shorter or simpler. In our
opinion, none of them reached the clarity of exposition in
\cite{5}.

Two years later, in 1958, Dyer and Hamstr\"{o}m applied this theorem to
get the sufficient conditions for a regular map $f$ to be a trivial fibration.
Such a condition turned out to be local $n$-connectedness ($LC^{n}$) of the
homeomorphisms group $H(M)$ of the fiber $M$ of $f$. The problem
when $H(M)$ is $LC^{n}$ was one of the central in topology over a period of
almost 20 years and served as one of the key sources for the
development of
infinite--dimensional topology, as a separate part of topology.

For a first encounter with the theory of selections, the
papers
\cite{4} and \cite{5} are too difficult and too voluminous. On the
other hand, the short note \cite{3} quickly tells the reader of
the most popular method of selection theory -- the method of 
outside approximation. The note consists of the proof of the
Convex--valued and the the $0$--dimensional selection Theorems.
The last theorem is a particular case of the Finite-dimensional
theorem for $n=0$.

\proclaim{Theorem 4} If $X$ is zero-dimensional  ($\dim X=0$) and
paracompact space, and if $Y$ is a complete metric space, then every LSC
mapping $\phi:X \to \FF(Y)$ admits a singlevalued continuous
selection.
\endproclaim

In spite of its relative simplicity and clarity of its proof, the
Zero-dimensional selection theorem has many applications in selection
theory and other areas of mathematics.

The last paper of the series \cite{4,5,6}  dealt mainly
with  restrictions on the displacement of a closed subset $A$ in
$X$. For example, as in the Borsuk pairs, when $X=Z \times [0;1]$
and $A=(Z \times {0})\cup (B \times [0;1])$ for an appropriate $B
\subset X$. Also the lower semicontinuity assumption in \cite{6}
was strengthened by continuity in the
corresponding Hausdorff metric $h_{\rho}$ in exp $X$.
Here we reproduce a typical statement, Theorem 6.1:

\proclaim{Theorem 5} Let $X$ be a paracompact space with $\dim X
\leq n+1$ and $A\subset X$  a weak deformation retract of $X$. Let
$(Y,\rho)$ be a complete metric space, $\FF$  a uniformly-$LC^{n}$
family of nonempty closed subsets of $(Y,\rho)$ and $\phi:X \to
\FF$ a continuous map with respect to $h_{\rho}$. Then every selection
of $\phi|_{A}$ can be extended to a singlevalued continuous
selection of $\phi$.
\endproclaim

Surprisingly deep
constructions and results of \cite{6}  have until now had no real
and clear applications.

\head {\bf 4. Papers from 1959} \endhead

We begin by the first paper of the series \cite{7,8,9,10}. If one
combines arbitrary paracompact domains, as in the Convex-valued
selection theorem, and arbitrary complete metric ranges for
closed-valued mappings, as in the Zero-dimensional selection
theorem, then of course, there is no hope of obtaining a {\it
singlevalued} continuous selection. It turned out that under those
assumptions a sufficiently fine {\it multivalued} selections exist.
It was rather an unexpected and "...curious result about
semi-continuous..., \cite{7}" selections. Below, $2^{Y}$ denotes
the family of all nonempty subsets of a set $Y$:

\proclaim{Theorem 6 (\cite{7; Theorem 1.1})} Let $X$ be a paracompact space, $Y$
a metric space, and $\phi:X \to 2^{Y}$ an LSC map with each $\phi(x)$
complete. Then there exist $\psi:X \to 2^{Y}$ and $\theta:X \to
2^{Y}$ such that:
\itemitem{(a)} $\psi(x) \ss \theta(x) \ss \phi(x)$ for all $x \in X$;
\itemitem{(b)} $\psi(x)$ and $\theta(x)$ are compact, for all $x \in
X$;
\itemitem{(c)} $\psi$ is LSC; and
\itemitem{(d)} $\theta$ is USC (upper semicontinuous).
\endproclaim

It appears
that this was in principle, the very
first theorem on multivalued
selections. The proof of this Compact-valued selection theorem is
based on the so called method of 
inner approximations. Roughly
speaking, one can inscribe
into each value $\phi(x)$
a tree with a
countable set of levels, with finite sets of vertices
on each
level so that each maximal linearly ordered sequence of vertices
will be fundamental.

Thus the sets $\psi(x)$ and $\theta(x)$ are
constructed as the sets of limits of different kinds of such
maximal paths in the tree. Shortly,
$\psi(x)$ and $\theta(x)$ are
limits of certain
inverse (countable) spectra in the complete metric
space $\phi(x)$.
Beginning by \cite{7} multivalued selections became
by then a fully respected part
of general selection theory.
The comprehensive  fundamental paper \cite{10} also had an important
impact on the development of selection theory. 

In this paper
the axiomatic theory of convexity in metric spaces was presented.
As far as we know,
this was also one of the first papers on
axiomatic convexities. It served as the
starting point for
many investigations in this direction.

Also, the method of inner approximations from \cite{7} was changed
and applied in \cite{10} to convex-valued maps. Roughly speaking,
at each level of a tree above one can consider the barycenter of
all vertices at that level, with respect to a suitable continuous
partition of the unity of the domain. In this way it is possible
to obtain a poinwise convergent sequence of singlevalued
(discontinuous!) selections with degree of discontinuity uniformly
tending to zero. Therefore the limit gives the desired continuous singlevalued
selection.

In our experience, we have encountered several times the situations
when the simpler and more direct smoothing method of  outside
approximations did not work, whereas the method of  inner
approximations successfully solved the problem at hand. Looking at the data  on submission of the
papers, one may perhaps infer that \cite{10} was originally the source for \cite{7}.

Whereas Lemmas 5.1 and 5.2 and Theorem $3.1'''$ were proved in
\cite{4} for perfectly normal domains and separable Banach range
spaces,   a  version was obtained in
\cite{8} for metric domains and any Banach range spaces. The proof
was based on the replacement of the $G_{\delta}$-property for
closed subsets of a perfectly normal domain by the A. Stone
theorem on the existence of $\sigma$-discrete closed basis in any
metric space. Note also that Theorem 5.1 \cite{8} on the one hand,
used the ideas from the proofs in \cite{6}, and on the other hand
was the basis for the later appearance of such notions as SEP and
SNEP ({\it selection extension } and {\it selection neighborhood
extension properties}) \cite{18}.

While \cite{10} estimates the relations and links between convex and
metric structures on the set, the paper \cite{9} deals with the
degree of nonconvexity of a closed subset $P$ of a Banach space,
endowed with standard convex and metric structures. Simply put,
imagine that we move the
endpoints of a segment of length $2r$ over a
set $P$. In this situation it is very natural to look for the
distance between the points of segment and the set $P$.

So if all such distances are less than or equal to $\a \cdot r$
for some constant $\a \in [0;1)$, then the set $P$ is {\it
paraconvex} in dimension $1$. By passing to triangles,
tetrahedra, and $n$-simplices, one
obtains the notion of a
{\it
paraconvex}
set. So, as was proved in \cite{9},
the statement
of the Convex-valued selection theorem \cite{3,4} holds whenever
one replaces the convexity assumption for the values $\phi(x)$ by
their $\a$-paraconvexity, for some common $\a \in [0;1)$, for all
$x \in X$.

Moreover, the proof looks as a double sequential ``improvement''
process of exactness of approximation, on the account of applying
the Convex-valued selection theorem.  

\head {\bf 5. Papers from 1964--1979 } \endhead

One of the main purposes of the series \cite{12-15} was to examine
the metrizability assumption for the range space in the
Convex-valued selection theorem. In the papers  \cite{12,13,14}
improvements of the Arens-Eells embedding theorem
were proved and a selection theorem for mappings from
metric domains into completely metrizable subsets of locally
convex topological vector (LCTV) spaces  was established.
It was shown in
\cite{15} that the statement holds for paracompact domains as
well. Observe that for LCTV spaces {\it completness} is a delicate and
in general, "multivalued" notion. Below, a LCTV space is said to be {\it
complete} if the closed convex hull of any compact subset is also
a compact subset.

\proclaim{Theorem 7 (\cite{15; Theorem 1.2})} Let $X$ be a
paracompact space
and $(M,\rho)$ a metric subset of a complete LCTV space $E$. Let $\phi:
X \to 2^{M}$ be an LSC map
such that every $\phi(x)$ is $\rho$-complete.
Then there exists a continuous singlevalued $f:X \to E$ such that
for every $x \in X$, the value $f(x)$ belongs to the closed convex
hull of the set $\phi(x)$.
\endproclaim

Note that one of the key ingredients of the proof is the
Compact-valued selection theorem.
Next, if $\phi$ is
convex-valued and closed-valued, then completness of the entire
$E$ can be replaced by completness of the closed spans of
$\phi(x), x \in X$. Such a replacement can be also derived
from the Zero-dimensional selection theorem and by the technique
of pointwise integration (see  \cite{R13}).  

In the joint paper  with Engelking and Heath \cite{16},  Michael in
some sense returned to his first selection publication \cite{1}.
Namely, by using embedings into closed topologically
well-ordered subspaces of the Baire space $B({\goth m})$, they
proved (\cite{16; Corollary 2}) that for any complete metric,
zero-dimensional (with respect to $\dim$ or  Ind) space
$(X,\rho)$ there exists a singlevalued continuous {\it selector}
$f$ on the family $\FF(X)$ of all nonempty closed subsets of $X$.

Here $\FF(X)$ is endowed with the Hausdorff topology, say
$\tau_{\rho}$, and $f:\FF(X) \to X$ is a mapping with $f(A) \in A$
for every $A \in \FF(X)$.  The zero-dimensionality is the
necessary restriction, because there are no selectors for
$\FF(\R)$ (see  \cite{16; Proposition 5.1}).

Note that formally, a selector is simply a selection of the
multivalued evaluation mapping, which associates to each $A \in
\FF(X)$ the same $A$,  but as a subset of $X$. However,
historically the situation was reverse. In \cite{1} Michael
proposed a separation of the problem about the
existence of a
selection $g: Y \to X$ for $G:Y \to 2^{X}$
into two separate problems: first, to check
that $G$ is continuous and second, to prove that there exists a
selector on $2^X$. Hence, the selection problem was originally
reduced to a certain selector problem.  

\head {\bf 6. Papers from 1980--1990} \endhead

Pick points $x_{1}, x_{2},...,x_{n}$  in the
domain $X$ of a
multivalued mapping $\phi$ and arbitrary select points $f(x_{i})
\in \phi(x_{i})$, using the
 Axiom of choice. Thus we find a
partial selection of $\phi$ over the
closed subset $C=\{x_{1},
x_{2},...,x_{n}\} \ss X$. By replacing the values $\phi(x_{i})$
with
the singletons $\{f(x_{i})\}$ we once again obtain an LSC mapping,
say $\phi_C$. If all assumptions of a  selection theorem
hold for the new LSC mapping $\phi_C$, then such a mapping admits
a selection, and hence $\phi$ also admits a selection.

This simple observation shows that any restriction
for the value
of $\phi$ over a finite subset $C \ss X$, like closedness,
connectivity, convexity, etc. are inessential for the existence of
a continuous selection of $\phi$. But what can one say about such
an omission for an
{\it infinite} $C \ss X$? Clearly, $C$ should be a
sufficiently "small", "dispersed", etc. subset of $X$. At the
international congress of mathematicians in  Vancouver in 1974, Michael
announced results for  countable $C$. Based on this, the following
result was published in 1981 (see \cite{19; Theorem 1.4}):

\proclaim{Theorem 8} Let $X$ be a paracompact space, $Y$ a Banach
space, $C \ss X$ a countable subset and $\phi: X \to 2^{Y}$ an LSC
map with closed and convex values $\phi(x)$ for all $x \notin C$.
Then for every closed subset $A \ss X$, each selection of
$\phi|_{A}$ admits an extension which is a selection of $\phi$
(shortly, $\phi$ has SEP).
\endproclaim

Briefly, over a countable subset of a domain we can simply omit
any restriction for the
values of LSC $\phi: X \to 2^{Y}$. A year
before, in a joint paper
with Pixley  \cite{17}, Michael proved that the
convexity assumption can be omitted over any subset $Z \ss X$ with
$\dim_{X}Z=0$.

Roughly speaking, results of \cite{17,18,19,23,25} are principally
related to  several possibilities for relaxing convexity
in  selection
theorems and in particular, the  closedness assumptions
for values of multivalued mappings. For example, let us mention
the following two results:

\proclaim{Theorem 9 (\cite{19; Theorem 7.1})} Let $X$ be a
paracompact space, $Y$ a Banach space ,$C \ss X$ a countable
subset, $Z \ss X$ a subset with $\dim_{X}Z \leq 0$ and $\phi: X \to
2^{Y}$ an LSC map such that $\phi(x)$ is closed for all $x \notin
C$ and $Clos(\phi(x))$ is convex, for all $x \notin Z$. Then
$\phi$ has SEP.
\endproclaim

\proclaim{Theorem 10 (\cite{18; Theorem 1.2})} Let $X$ be a
paracompact space, $Y$ a Banach space,  $Z \ss X$ a subset with
$\dim_{X}Z \leq n+1$ and $\phi: X \to \FF(Y)$ an LSC map such that
and $\phi(x)$ is convex, for all $x \notin Z$ and the family
$\{\phi(x): x \in Z\}$ is uniformly equi-$LC^{n}$. Then $\phi$ has
SNEP. If moreover, $\phi(x)$ is $n$-connected for every $x \in Z$,
then $\phi$ has SEP.
\endproclaim

Note that in \cite{18} the technique of the proof in \cite{5}
was rearranged in a more structured form, with exact extracting
of the
useful properties like SEP, SNEP and SAP ({\it selection
approximation property}).

The joint paper with Magerl and Mauldin  \cite{20} formally
contains no "selections" in the title or
in the statements of the
main theorems (1.1 and 1.2). Nevertheless, the essence of these
theorems is contained in the selection result.

It is a classical
fact that each metric compact $X$ can be represented as the
image of
the Cantor set $K$ under some continuous surjection $h: K \to X$.
Theorem 5.1 of
\cite{20} states that if $\{X_{\a}\}$ is a family of
subcompacta of a metric space $X$ which is continuously
parameterized by $\a \in A$ with $\dim A=0$ then one can choose a
family of surjections $h_{\a}: K \to X_{\a}$ continuously
depending on the same parameter $\a \in A$. Such parameterized
version of the
Alexandrov theorem is in fact, derived from
the Zero-dimensional selection theorem.

In general, the decade 1980--1990 was marked by
Michael's very diverse set of papers on selections, practically
every one of which contained new ideas of high quality. For one
more example,  the Finite-dimensional selection
theorem from \cite{5} was strengthened in \cite{23}
simultaneously in two
directions.
First, the assumption that $\{\phi(x)\}_{x \in X}$
is an
equi-$LC^{n}$ family in $Y$ was
replaced by the property that
fibers  $\{\{x\} \times \{\phi(x)\}_{x \in X}\}$ constitutes an
equi-$LC^{n}$ family in $X \times Y$. This answered the problem
of Eilenberg stated in 1956 (see the
comments in \cite{5}).
Next, the
closedness assumption for $\phi(x) \ss Y$ can be weakened to the
closedness of graph-fibers $\{x\}\times \phi(x)$ in some
$G_{\delta}$-subset of $X \times Y$.

The key ingredient of the proof was a "factorization"
construction. Briefly, it occurred that the
LSC mapping $\phi: X \to
Y$ with weakened assumptions can be represented as a composition
$\phi=h \circ \psi$ with singlevalued $h:Z \to Y$ and with $\psi:X
\to Z$ satisfying 
the classical assumptions of the
Finite-dimensional selection theorem \cite{5}. Hence the
composition of a selection of $\psi$ with $h$ gives the desired
selection of $\phi$.

We guess that the idea of the
appearance of
the
$G_{\delta}$-conditions
was a corollary of constructions of selections, avoiding a
countable set of obstructions, from the paper
\cite{22} which appeared one year earlier:

\proclaim{Theorem 11 (\cite{22; Theorem 3.3})} Let $X$ be a
paracompact space,
$Y$ a Banach space and $\phi: X \to \FF(Y)$ a LSC map with convex
values. Let $\psi_{i}: X \to \FF(Y), i \in \N$ be continuous,
$Z_{i} = \{x \in X: \phi(x) \cap \psi_{i}(x) \not= \emptyset\}$
and suppose that
$$
dim\,X<dim\,\phi(x)-dim(conv(\phi(x) \cap \psi_{i}(x))),
$$
for all $x \in Z_{i}$ and
$i \in \N$. Then $\phi$ admits a
selection $f$ which
avoids every $\psi_{i}$: $f(x) \notin \psi_{i}(x)$.
\endproclaim

Briefly, in the values $\phi(x)$ there is sufficient "room"
to avoid all sets $\psi_{i}(x)$.

Based on \cite{22,23}, Michael
stated
in "Open problems in topology, I" the "$G_{\delta}$"-problem \cite{24; Problem 396}:
Does
the Convex-valued selection theorem remain true if $\phi$ maps $X$
into some $G_{\delta}$-subset $Y$ of a Banach space $B$ with
convex values which are closed in $Y$? In spite of numerous cases
with affirmative answer this problem has in general a negative (as
it was
conjected in \cite{24}) solution, for details see the paper of
Namioka and Michael in this issue.

\head {\bf 7. Papers from 1992} \endhead

In general, all papers \cite{21,25,27} are related to "dispersed",
mainly to zero-dimensional, (in $\dim$-sense) domains of
multivalued mappings.

Briefly, in \cite{25} results of \cite{17; Theorem 1.1} and \cite{19; Theorem 1.3}
are unified and generalized in the spirit of \cite{23} to  subsets
$C \ss X$, which are unions of countable family of
$G_{\delta}$-subsets $C_{n}$ of $X$ and to a mappings $\phi$,
having SNEP at each $C_{n}$. In the paper written with Choban
\cite{27},
the Compact-valued selection Theorem 6 was derived from
the Zero-dimensional one
(Theorem 4).
In fact, a paracompact domain
$X$ was presented as the
image $h(Z)$ of some zero-dimensional
paracompact space $Z$ with respect to some appropriate
continuous
(perfect or inductively open) mapping $h: Z \to X$.
Theorem 4 applied
to the composition $\phi\circ h$ gives a
selection, say $s:Z \to Y$. So, the composition $s \circ h^{-1}$
will be a desired multivalued selection of $\phi: X \to \FF(Y)$.

The pair of papers \cite{26,28} is related to "the differences
between selection theorems which assume that the domain is
finite-dimensional and those which do not". More generally, based
on the Pixley counterexample in \cite{26} it was shown that a
genuine dimension-free analogue of the Finite dimensional selection
theorem does not exist or  briefly, that there are no purely
topological analogs of convexity. In comparison, in \cite{28} a
convexity, or connectivity, type restrictions in the spirit of
\cite{10} for a {\it mapping} are presented and under such
restrictions
 the equivalence is proved
 between the existence of
global selections and the existence of
selections locally.

The paper \cite{29} on continuum-valued selections is an elegant
simultaneous application
of the
"universality" idea from \cite{27} and the
one-dimensional selection theorem (special case $n=0$ of
Theorem 5). The key step can be described as follows.
Due to a recent theorem of Pasynkov, each paracompact domain $X$ can
be represented in the form $h(Z)$, for some perfect, open surjection
$h:Z \to X$ with pathwise connected fibers and for some
paracompact space $Z$ with $\dim\,Z\leq 1$. So, if the composition
$\phi\circ h$ admits a selection, say $s:Z \to Y$ then the
composition $s \circ h^{-1}$ will be a continuum-valued selection
of $\phi: X \to \FF(Y)$.

We should mention the thoughtfulness, exactness, perfectness of
all Michael's papers. His laconic style of exposition is perfectly
matched with the
deepness of his results. In our opinion,
{\sl A man of few words but with great ideas} could well serve as a good description of his character.
As a rule, all his papers are equipped with a considerable
number of additional references, which were added at proofs, and
which very precisely give correct accents to the paper needed for
proper understanding.
In conclusion of this survey of Michael's results on
selections we wish our jubilant successful realization of many more projects.

\head {\bf Acknowledgements} \endhead
We thank Jan van Mill for comments and suggestions. The first author was supported by the Slovenian Research Agency
grants No. P1-0292-0101-04 and Bl-RU/05-07-04. The second author
was supported by the RFBR grant No.05-01-00993.

\enddocument